\documentclass[10pt,twocolumn,twoside]{IEEEtran}
\usepackage{amsmath,amsfonts,amssymb}
\usepackage{balance}
\usepackage{xcolor}
\usepackage{algpseudocode}
\usepackage{subfig}
\usepackage{cite}
\usepackage[pdftex]{graphicx}
\DeclareGraphicsExtensions{.jpg,pdf,.png}
\usepackage{amscd,algorithm,algorithmicx}
\ifCLASSINFOpdf
\else
\fi
\date{}
\begin{document}
% #################  Title  ######################################
\title{Recovery of Sparse and Low Rank Components of Matrices Using Iterative Method with Adaptive Thresholding}
% #################  Authors  ####################################
\author{Nematollah Zarmehi, \IEEEmembership{Student Member, IEEE} and Farokh Marvasti, \IEEEmembership{Senior Member, IEEE}\thanks{N. Zarmehi and F. Marvasti are with the Advanced Communication Research Institute (ACRI), Electrical Engineering Department, Sharif University of Technology, Tehran, Iran. (e-mail: zarmehi\_n@ee.sharif.edu, marvasti@sharif.edu)}}
% #################  Header  #####################################
%\markboth{}{Zarmehi \MakeLowercase{\text	\it{et al.}}: Recovery of Sparse and Low Rank Components of Matrices Using Iterative Method with Adaptive Thresholding}
\maketitle
% #################  Abstract  ###################################
\begin{abstract}
In this letter, we propose an algorithm for recovery of sparse and low rank components of matrices using an iterative method with adaptive thresholding. In each iteration, the low rank and sparse components are obtained using a thresholding operator. This algorithm is fast and can be implemented easily. We compare it with one of the most common fast methods in which the rank and sparsity are approximated by $\ell_1$ norm. We also apply it to some real applications where the noise is not so sparse. The simulation results show that it has a suitable performance with low run-time.
\end{abstract}

% #################  Index Terms  ################################

% #################  Introduction  ###############################
\section{Introduction}\label{sec:intro}
\IEEEPARstart{I}{n} some applications, we need to decompose a given matrix $\mathbf{Y}\in \mathbb{R}^{m\times n}$ into a low rank matrix $\mathbf{L_*}\in \mathbb{R}^{m\times n}$ and a sparse matrix $\mathbf{E_*}\in \mathbb{R}^{m\times n}$, i.e.,
\begin{equation}\label{eq:base}
\mathbf{Y} = \mathbf{L_*}+\mathbf{E_*},
\end{equation}
where the rank of $rank(\mathbf{L_*})=r_0 \ll \min(m,n)$ and $sparsity(\mathbf{E_*})=k_0 \ll mn$. This problem is arises in various applications such as subspace clustering \cite{ref:subsclust}, matrix rigidity in computation complexity \cite{ref:matrixrigid}, removing shadows and specularities from sequence of images \cite{ref:face2,ref:canPCA,ref:face1}, background modeling in video surveillance \cite{ref:bg1}, link prediction in social networks \cite{ref:linkpred}, and graph clustering \cite{ref:graphclust}.

When the data is correlated and its dimension is high, it lies on a lower dimension linear space. This property is the base of some data analyses such as Principal Component Analysis (PCA) \cite{ref:pcaspringer}, independent component analysis \cite{ref:icabook}, and sparse component analysis \cite{ref:sparseca,ref:sparseca2}. For example, when the noise components (the entries of $\mathbf{E_*}$) in (\ref{eq:base}) are i.i.d. Gaussian, the PCA optimally decomposes the low rank and noise components. Speaking mathematically, in the case of Gaussian noise, the PCA gets the optimal solution by minimizing the norm of $\mathbf{E_*}\in \mathbb{R}^{m\times n}$ subject to a condition on the rank of $\mathbf{L_*}$ as follows \cite{ref:pcaspringer}:
\begin{equation}\label{eq:classicalPCA}
\begin{split}
\arg\min_\mathbf{E_*} & ~ \|\mathbf{E_*}\|, \\
s.t.: & ~ rank(\mathbf{L_*}) \leq r ~ \mbox{and} \\
&~ \mathbf{Y} = \mathbf{L_*} + \mathbf{E_*}.
\end{split}
\end{equation}
This optimization problem could be easily done by computing the Singular Value Decomposition (SVD) of $\mathbf{Y}$ and zeroing the $\min(m,n)-r$ right hand side Singular Values (SVs). Although this is very simple, it is well known that PCA may fail to do such decomposition when the noise components do not obey an i.i.d. Gaussian distribution. For example, if only one of the entries of low rank component faces with a large corruption, PCA fails to do decomposition successfully. These kind of corruptions are very usual in digital signal processing.

To overcome this issue, many efforts have been made to make the PCA robust against non-Gaussian noises. Among them, Pursuit Component Analysis (PCP) is an idealized method that guarantees the exact recovery of sparse and low rank components under some conditions. For this purpose, PCP solves the following convex optimization problem \cite{ref:canPCA}:
\begin{equation}\label{eq:canpca}
\begin{split}
\arg\min_{\mathbf{L_*},\mathbf{E_*}} & ~ \|\mathbf{L_*}\|_* + \lambda \|\mathbf{E_*}\|_1, \\
s.t.: &~ \mathbf{Y}=\mathbf{L_*}+\mathbf{E_*},
\end{split}
\end{equation}
where $\lambda >0$ is a weighting factor that balances the sparsity and low rank. Here, the rank and $\ell_0$-norm are approximated by the nuclear norm ($\|.\|_*$) and $\ell_1$-norm, respectively, to make the primary non-convex optimization problem convex.

In \cite{ref:tesnor}, an iterative hard thresholding method is proposed for tensor recovery under the affine constraint that is different from the mathematical model introduced in (\ref{eq:base}). It is a variant of normalized iterative hard thresholding proposed in \cite{ref:niht} for sparse recovery. An iterative hard thresholding is also proposed for matrix completion in \cite{ref:mtxccom}. This method also obeys an affine condition.

In this letter, we consider the model introduced in (\ref{eq:base}) and also suppose that the noise components are not necessarily i.i.d. Gaussian and we don't know the rank of $\mathbf{L_*}$ and sparsity of $\mathbf{E_*}$. For example, $\mathbf{E_*}$ can be sparse. We propose an iterative method with adaptive thresholding to recover the sparse and low rank components. The idea of adaptive thresholding is adopted from the proposed method in \cite{ref:imat} for recovery of the original signal from its random samples.

The rest of this letter is organized as follows: Section \ref{sec:alg} is devoted to propose our method for recovery of sparse and low rank components of matrices. In Section \ref{sec:sim}, the simulation results are presented and the proposed method is evaluated with different scenarios. We also apply our proposed method on some real applications in which the sparse and low rank recovery is needed. At last, Section \ref{sec:conclusion} concludes the letter.

\section{The Proposed Algorithm}\label{sec:alg}
Considering the mathematical model (\ref{eq:base}), we aim to recover $\mathbf{L_*}$ and $\mathbf{E_*}$ from the given matrix $\mathbf{Y}$. Since the entries of $\mathbf{E_*}$ and SVs of $\mathbf{L_*}$ are sparse, we propose to threshold them adaptively during iterations until the stop condition is satisfied. The adaptive threshold level of $k$-th iteration is suggested as follows:
\begin{equation}
\tau_k = \beta \sigma_1 e^{-\alpha k},
\end{equation}
where $\alpha$ and $\beta$ are two constants and $\sigma_1$ is the largest SV of $\mathbf{Y}$. On each iteration, we force the low rank and sparse terms to satisfy (\ref{eq:base}). Then, our algorithm has the following steps:
\begin{equation}
\left\{ \begin{array}{l}
\mathbf{L} = \mathcal{T}_{\tau_k}(\mathbf{L}_{i-1})\\
\mathbf{E}=\mathbf{Y}-\mathbf{L}\\
\mathbf{E}_i=\mathcal{D}_{\tau_k}(\mathbf{E})\\
\mathbf{L}_i=\mathbf{Y}-\mathbf{E}_i.
\end{array} \right.
\end{equation}
$\mathcal{T}_{\tau_k}(\mathbf{X})$ with $\mathbf{X}=\mathbf{U}\mathbf{\Sigma}\mathbf{V}^T, \mathbf{X}\in \mathbb{R}^{m\times n}$ and $\mathcal{D}_{\tau_k}(\mathbf{X})$ are defined as follow:
\begin{equation}
\begin{split}
\mathcal{T}_{\tau_k}(\mathbf{X}) & \buildrel \Delta \over = \mathbf{U}\mathcal{D}_{\tau_k}(\mathbf{\Sigma})\mathbf{V}^T\\
\mathcal{D}_{\tau_k}(\mathbf{X}) & \buildrel \Delta \over = \left[{\{x_{i,j}-\tau_k\}^+}\right]_{i=1,j=1}^{m,n},
\end{split}
\end{equation}
where $\{x\}^+=\max(x,0)$. Here, a matrix $\mathbf{X}\in \mathbb{R}^{m\times n}$ is denoted by $[x_{i,j}]_{i=1,j=1}^{m,n}$.

We name our method Sparse and Low rank Recovery using Iterative Method with Adaptive Thresholding (SLR-IMAT) that is shown in Algorithm \ref{alg1}. The proposed method consists of two loops. At $k$-th iteration of the outer loop, the inner loop iterates $M2$ times with a threshold level of $\tau_k$. The outer loop is terminated until the stop criterion is satisfied. We terminate the it when $\|\mathbf{\hat{L}_{k}}-\mathbf{\hat{L}_{k-1}}\|_F<\epsilon$ or it reaches the maximum number of iterations. Here, $\mathbf{\hat{L}}_k$ is the estimated low rank matrix at $k$-th iteration.

\begin{algorithm}
\caption{SLR-IMAT}\label{alg1}
\begin{algorithmic}[1]
\State \textbf{input:}
\State \indent {Data matrix:} $\mathbf{Y} \in \mathbb{R}^{m\times n}$
\State \indent {Stopping threshold:} $\epsilon$
\State \indent {Two constants:} $\alpha, \beta$
\State \indent {Max. number of iterations of the inner loops:} $M1$
\State \indent {Max. number of iterations of the outer loops:} $M2$
\State \textbf{initialization:}
\State \indent $k \gets 0$
\State \indent $\mathbf{L} \gets \mathbf{Y}$
\State \indent $\mathbf{\hat{L}_0}=\mathbf{Y}$
\While {$e>\epsilon ~ \& ~k<M1$}
\State $\tau_k \gets \beta \sigma_1(\mathbf{Y})e^{-\alpha k}$
\For {$i=1:M2$}
\State $[\mathbf{U},\mathbf{\Sigma},\mathbf{V}]=svd(\mathbf{L})$
\State $\mathbf{\Sigma}(\mathbf{\Sigma}<\tau_k)=0$
\State $\mathbf{E} = \mathbf{Y} - \mathbf{U}\mathbf{\Sigma}\mathbf{V^T}$
\State $\mathbf{E}(|\mathbf{E}|<\tau_k)=0$
\State $\mathbf{L} \gets \mathbf{Y}-\mathbf{E}$
\EndFor
\State $k\gets k+1$
\State $\mathbf{\hat{L}_k}=\mathbf{L}$
\State $e=\|\mathbf{\hat{L}_{k}}-\mathbf{\hat{L}_{k-1}}\|_F$
\EndWhile
\State \indent $\mathbf{L_*}\gets \mathbf{\hat{L}_k}$
\State \indent $\mathbf{E_*}\gets \mathbf{Y}- \mathbf{\hat{L}_k}$
\State \textbf{return} $\mathbf{L_*}, \mathbf{E_*}$
\end{algorithmic}
\end{algorithm}

\section{Simulation Results}\label{sec:sim}
This section presents the simulation results. The simulations are done by the MATLAB R2015a software on Intel(R) Core(TM) i7-5960X @ 3GHz system with 32GB-RAM. First, we shall compare the proposed method with the Inexact Augmented Lagrange Multiplier (IALM) method introduced in \cite{ref:inexact} as a method that solves the minimization problem suggested in \cite{ref:canPCA}. This method is a fast variant of the Augmented Lagrange Multiplier (ALM) for sparse and low rank recovery introduced in \cite{ref:alm1}. Then, the phase transition plot of the proposed method will be shown in two different scenarios. Finally, we apply the proposed method to two real applications in which the sparse and low rank recovery is used. We use $SNR(\mathbf{L_*},\mathbf{\hat L}) = 20\log_{10}\left(\|\mathbf{L_*}\|_F/\|\mathbf{L_*}-\mathbf{{\hat L}}\|_F\right)$ in dB to evaluate the performance of the recovery task. 

\subsection{Exact Recovery}\label{subsec:ex}
In this subsection, we verify that the proposed method can exactly separate the low rank and sparse components of a given matrix. For this purpose, at first, $r$-rank square matrices of size $m=n=500,1000,\cdots,3000$ are produced as the product of two matrices like $\mathbf{L_*}=\mathbf{A}\mathbf{B}^T$ where both $\mathbf{A}$ and $\mathbf{B}$ are independently sampled from a $\mathcal{N}(0,1/n)$ distribution. Then, the sparse noise component is added to the low rank component additively. The support set of sparse noise $\mathbf{E_*}$ is generated uniformly at random whose values have random signs. In other words, the non-zero values are independently selected from the set $\{-1,+1\}$ with equal probabilities. Two scenarios are considered in this subsection. In the first scenario, we set $k_0=0.05n^2$ which is supposed to be easier than the second scenario in which $k_0$ is set to $0.1n^2$. We compared the proposed method with the IALM method \cite{ref:inexact}. The IALM MATLAB code is downloaded from \textit{http://perception.csl.illinois.edu/matrix-rank/sample\_code.html} and the all parameters are set to their default values. The results are reported in Tables \ref{tab:ex1} and \ref{tab:ex2}.

\begin{table}[t!]
\centering
\caption{Exact recovery of low rank and noise components for random problems. The noise components are independent of the low rank components and have random signs. In this scenario, we set $r_0=0.05\times n$ and $k_0=0.05\times n^2$.}\label{tab:ex1}
\renewcommand{\arraystretch}{1.1}
\begin{tabular}{|c||c|c|c|c|c|} \hline
\textbf{Dimension} & $\mathbf{SNR_{i}}$ & \multicolumn{2}{|c|}{$\mathbf{SNR_{o} (dB)}$} & \multicolumn{2}{|c|}{\textbf{Time (s)}} \\  \cline{3-6}
$m = n$ & \textbf{(dB)} & \textbf{SLR-IMAT} & \textbf{IALM} & \textbf{SLR-IMAT} & \textbf{IALM} \\ \hline \hline
500 & -27.1 & 299.8 & 243.5 & 2.0 & 11.5  \\ \hline
1000 & -30.0 & 297.1 & 244.5 & 9.9 & 27.7 \\ \hline
1500 & -31.7 & 295.3 & 245.4 & 26.9 & 60.1 \\ \hline
2000 & -33.0 & 294.1 & 244.1 & 62.3 & 104.4 \\ \hline
2500 & -34.0 & 293.2 & 244.4 & 124.3 & 185.1 \\ \hline
3000 & -34.8 & 292.4 & 242.6 & 216.3 & 297.1 \\ \hline
\end{tabular}
\end{table}

In these Tables, $SNR_i$ indicates the input $SNR$, i.e., $SNR(\mathbf{L_*},\mathbf{E_*})$ while $SNR_o$ indicates the recovered $SNR$, i.e., $SNR(\mathbf{L_*},\mathbf{\hat{L}})$. According to the results of Tables \ref{tab:ex1} and \ref{tab:ex2}, both the proposed method and the IALM method could exactly recover the low rank and noise components but the run-time of the proposed method is lower than the run-time of the IALM method.

\begin{table}[t!]
\centering
\caption{Exact recovery of low rank and noise components for random problems. The noise components are independent of the low rank components and have random signs. In this scenario, we set $r_0=0.05\times n$ and $k_0=0.1\times n^2$.}\label{tab:ex2}
\renewcommand{\arraystretch}{1.1}
\begin{tabular}{|c||c|c|c|c|c|} \hline
\textbf{Dimension} & $\mathbf{SNR_{i}}$ & \multicolumn{2}{|c|}{$\mathbf{SNR_{o} (dB)}$} & \multicolumn{2}{|c|}{\textbf{Time (s)}} \\  \cline{3-6}
$m = n$ & \textbf{(dB)} & \textbf{SLR-IMAT} & \textbf{IALM} & \textbf{SLR-IMAT} & \textbf{IALM} \\ \hline \hline
500 & -31.7 & 271.2 & 245.6 & 2.0 & 10.2  \\ \hline
1000 & -34.7 & 284.3 & 248.5 & 9.8 & 25.0 \\ \hline
1500 & -36.5 & 287.8 & 247.7 & 26.4 & 54.0 \\ \hline
2000 & -37.8 & 287.9 & 243.3 & 62.0 & 101.1 \\ \hline
2500 & -38.7 & 287.5 & 240.0 & 124.1 & 183.5 \\ \hline
3000 & -39.5 & 286.9 & 243.7 & 216.6 & 302.3 \\ \hline
\end{tabular}
\end{table}

\subsection{Difficult Scenarios}\label{subsec:comp}
In this subsection, we consider more difficult scenarios than the scenarios investigated for the exact recovery. The low rank and noise components are produced in the same manner as explained in Subsection \ref{subsec:ex}. We set $k_0=0.3n^2$ and $0.4n^2$ in the first and in the second scenarios, respectively. It is obvious that these scenarios are more difficult than the scenarios of Subsection \ref{subsec:ex}. The results of our method and the IALM method are shown in Tables \ref{tab:comp1} and \ref{tab:comp2}. As it can be seen, our method has better performance than the IALM method because it could better recover the low rank and sparse components in terms of SNR even in lower time. We can see that the IALM method does not have suitable performance in the second scenario where $k_0=0.4n^2$. On the other hand, our method has about SNR greater than $100 dB$ in the second scenario.

\begin{table}[t!]
\centering
\caption{Comparison of the proposed method with the IALM method \cite{ref:inexact} for random problems. The noise components are independent of the low rank components and have random signs. In this scenario, we set $r_0=0.05\times n$ and $k_0=0.3\times n^2$.}\label{tab:comp1}
\renewcommand{\arraystretch}{1.1}
\begin{tabular}{|c||c|c|c|c|} \hline
\textbf{Dimension} & \multicolumn{2}{|c|}{\textbf{SNR (dB)}} & \multicolumn{2}{|c|}{\textbf{Time (s)}} \\  \cline{2-5}
$m = n$ & \textbf{SLR-IMAT} & \textbf{IALM} & \textbf{SLR-IMAT} & \textbf{IALM} \\ \hline \hline
500 & 127.1 & 66.5 & 2.7 & 34.1  \\ \hline
1000 & 147.5 & 81.0 & 35.0 & 107.0 \\ \hline
1500 & 150.4 & 79.4 & 68.5 & 250.2 \\ \hline
2000 & 151.1 & 82.6 & 138.3 & 447.2 \\ \hline
2500 & 152.0 & 84.0 & 222.0 & 692.1 \\ \hline
3000 & 152.5 & 82.3 & 322.9 & 986.1 \\ \hline
\end{tabular}
\end{table}

\begin{table}[t!]
\centering
\caption{Comparison of the proposed method with the IALM method \cite{ref:inexact} for random problems. The noise components are independent of the low rank components and have random signs. In this scenario, we set $r_0=0.05\times n$ and $k_0=0.4\times n^2$.}\label{tab:comp2}
\renewcommand{\arraystretch}{1.1}
\begin{tabular}{|c||c|c|c|c|} \hline
\textbf{Dimension} & \multicolumn{2}{|c|}{\textbf{SNR (dB)}} & \multicolumn{2}{|c|}{\textbf{Time (s)}} \\  \cline{2-5}
$m = n$ & \textbf{SLR-IMAT} & \textbf{IALM} & \textbf{SLR-IMAT} & \textbf{IALM} \\ \hline \hline
500 & 97.8 & 9.5 & 2.2 & 4.5  \\ \hline
1000 & 102.3 & 10.2 & 45.1 & 91.8 \\ \hline
1500 & 103.9 & 10.3 & 93.1 &  205.7  \\ \hline
2000 & 104.4 & 10.3 & 161.4 &  344.2  \\ \hline
2500 & 104.8 & 10.4 & 247.9 &  557.7  \\ \hline
3000 & 105.1 & 10.8 & 335.0 &  808.0 \\ \hline
\end{tabular}
\end{table}

\subsection{Phase Transition Plot}\label{subsec:phtr}
The phase transition plot is one of the most useful plots in low rank and sparse recovery problems. The phase transition plot shows the success recovery rate of the method for a wide range of rank and sparsity. Here, the low rank matrices are also generated as the same approach explained in Subsection \ref{subsec:ex}. The size of matrices is $200\times 200$ in this subsection. Again, we consider two scenarios. In the first scenario, noise is random while in the second, it is coherent with the low rank component. In random case, each noise component takes on values of $\{0,-1,+1\}$ with probabilities of $\{1-p,p/2,p/2\}$. In coherent case, the support set is the same as the random case but the values are the sign of the low rank components. To obtain the phase transition plot, we swap both $r_0/n$ and $p$ between 0.01 and 0.5 and each experiment is repeated 50 times. An experiment is declared successful if the recovered low rank component gets SNR greater than $60 dB$, i.e., $SNR(\mathbf{L_*},\mathbf{\hat L}) \geq 60 dB$.

\begin{figure}[!t]
\centering
\subfloat[Random Noise]{\includegraphics[width=.72\linewidth]{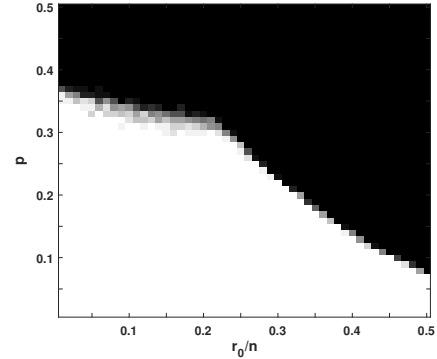}\label{subfig1:ph-tr-ex2}}
\hfil
\subfloat[Coherent Noise]{\includegraphics[width=.72\linewidth]{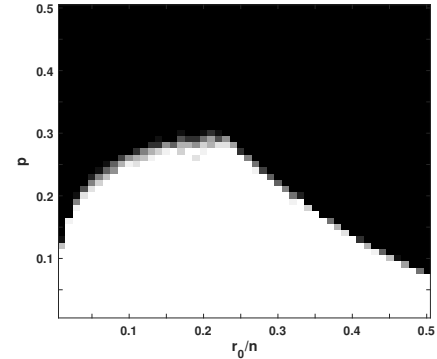}\label{subfig2:ph-tr-ex2}}
\caption{Phase transition between rank and sparsity. (a) In this case the noise is random and independent of the low rank component. (b) In this case the noise is coherent with the low rank component.}\label{fig:ph-tr-ex2}
\end{figure}

The success rate is shown by the gray color. The white and black colors indicate the success and failure, respectively. We can see the proposed method has acceptable performance. The proposed method has a vast region of success even in case of the coherent noise. 

\subsection{Applications}\label{sec:apps}
The separation of low rank and sparse components is a key approach in many applications. Here, we shall apply the proposed method on two of these applications, background modeling and removing shadows and specularities from face images. 

\begin{figure}[!t]
\centering
\subfloat{\includegraphics[width=.3\linewidth]{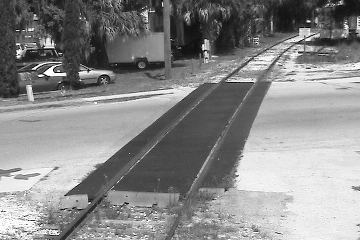}}
\hfil
\subfloat{\includegraphics[width=.3\linewidth]{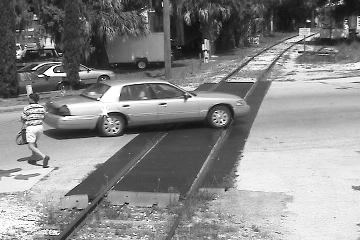}}
\hfil
\subfloat{\includegraphics[width=.3\linewidth]{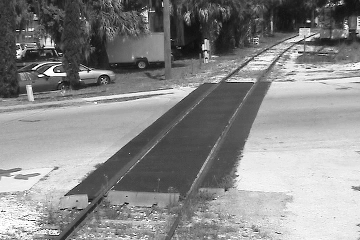}}\\
\subfloat{Frame \#1} \hfil \subfloat{~~~Frame \#2~~~} \hfil \subfloat{Frame \#3}\\
\vspace{0.2cm}
(a) Original frames \medskip \\ 
\subfloat{\includegraphics[width=.3\linewidth]{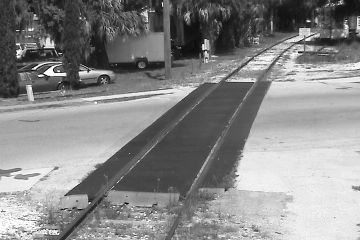}}
\hfil
\subfloat{\includegraphics[width=.3\linewidth]{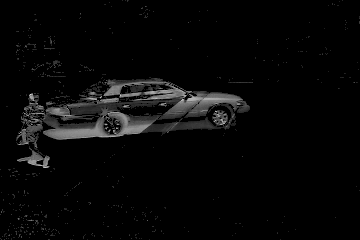}} \\
\subfloat{Background~~} \hfil \subfloat{~~Foreground}\\
\vspace{0.2cm}
(b) Results of SLR-IMAT \medskip \\
\subfloat{\includegraphics[width=.3\linewidth]{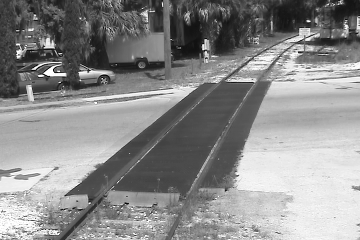}}
\hfil
\subfloat{\includegraphics[width=.3\linewidth]{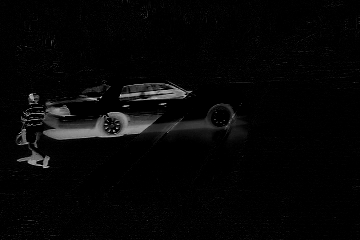}} \\
\subfloat{Background~~} \hfil \subfloat{~~Foreground}\\
\vspace{0.2cm}
(c) Results of IALM \medskip \\
\caption{An example of video sequence introduced in \cite{ref:nominalsequence}. (a) The original frames: Frames \#1 and \#3 are two sample frames from before and after the frame \#2 that contains the moving object. (b) The background and foreground obtained after low rank and sparse separation.}\label{fig:bg1}
\end{figure}

Background modeling is used in video surveillance applications. For example, consider a sequence of video frames from the same scene in which an object is moving. If we vectorize all frames and put them in a matrix, this matrix is low rank. The background is modeled as the low rank component and the foreground including the moving object is modeled as the sparse component. The proposed method and the IALM method are applied to two video sequences and the results are shown in Fig. \ref{fig:bg1}. In this figure, two sample frames from before and after the moving object and a sample frame containing the moving object are shown. It can be seen that the proposed method could successfully separate the background and foreground of two video sequences. 

\begin{figure}[!t]
\centering
\subfloat{\includegraphics[width=.17\linewidth]{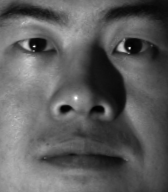}}
\hfil
\subfloat{\includegraphics[width=.17\linewidth]{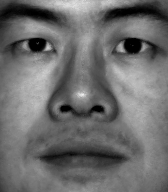}}
\hfil
\subfloat{\includegraphics[width=.17\linewidth]{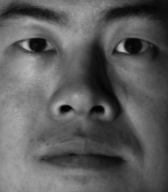}}
\\
\subfloat{\includegraphics[width=.17\linewidth]{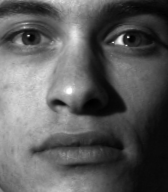}}
\hfil
\subfloat{\includegraphics[width=.17\linewidth]{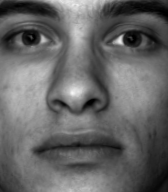}}
\hfil
\subfloat{\includegraphics[width=.17\linewidth]{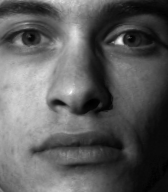}}
\\
\subfloat{\includegraphics[width=.17\linewidth]{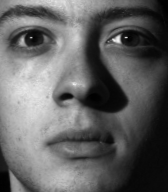}}
\hfil
\subfloat{\includegraphics[width=.17\linewidth]{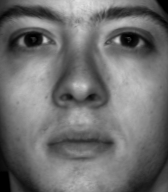}}
\hfil
\subfloat{\includegraphics[width=.17\linewidth]{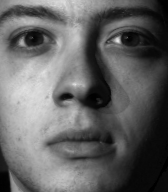}}
\\
\subfloat{\includegraphics[width=.17\linewidth]{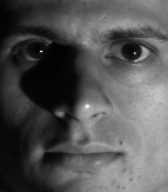}}
\hfil
\subfloat{\includegraphics[width=.17\linewidth]{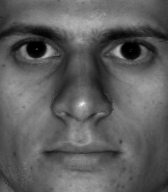}}
\hfil
\subfloat{\includegraphics[width=.17\linewidth]{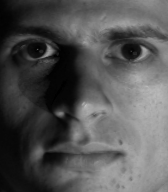}}
\\
\subfloat{(a)} \hfil \subfloat{~~~~(b)~~~~} \hfil \subfloat{(c)} \\
\caption{Removing shadows and specularities from the face images. (a) Original face images, (b) results of the SLR-IMAT method, and (c) results of the IALM method.}\label{fig:face}
\end{figure}

The low rank and sparse recovery are also used in removing shadows and specularities from face images \cite{ref:canPCA}. We have used the extended Yale face database B \cite{ref:exyale} for this simulation. Again, if we vectorize all face images and put them in a matrix the resultant matrix has low rank and the shadows and specularities and other noises would be canceled by low rank and sparse separation. Four images are selected from the database and the result of the proposed method and the IALM method are shown in Fig. \ref{fig:face}. Comparing the result images with the original ones, we can see that the shadows and specularities of the original faces are removed. Moreover, the output images of the IALM method is darker than the output images of the proposed method.

\section{Conclusion}\label{sec:conclusion}
We proposed a method for recovery of sparse and low rank components of matrices in this letter. The proposed method is iterative and uses an adaptive thresholding operator in each iteration to achieve low rank and sparse matrices. We compared the proposed method with the IALM method in which the $\ell_0$ norm is replaced by the $\ell_1$ norm. According to the simulation results, the proposed method outperforms the IALM method in terms of both SNR and run-time. Moreover, in the case of almost exact recovery ($k_0=0.05n, 0.1n$), the run-time of our method is lower than the run-time of IALM method. We have also applied the proposed method to two real applications that the interfering component was not so sparse and it could successfully separate the low rank and the interfering components.

% #################  References  #################################
\ifCLASSOPTIONcaptionsoff
  \newpage
\fi
\bibliographystyle{ieeetran}
\bibliography{Citations}

\begin{thebibliography}{10}
\providecommand{\url}[1]{#1}
\csname url@samestyle\endcsname
\providecommand{\newblock}{\relax}
\providecommand{\bibinfo}[2]{#2}
\providecommand{\BIBentrySTDinterwordspacing}{\spaceskip=0pt\relax}
\providecommand{\BIBentryALTinterwordstretchfactor}{4}
\providecommand{\BIBentryALTinterwordspacing}{\spaceskip=\fontdimen2\font plus
\BIBentryALTinterwordstretchfactor\fontdimen3\font minus
  \fontdimen4\font\relax}
\providecommand{\BIBforeignlanguage}[2]{{%
\expandafter\ifx\csname l@#1\endcsname\relax
\typeout{** WARNING: IEEEtran.bst: No hyphenation pattern has been}%
\typeout{** loaded for the language `#1'. Using the pattern for}%
\typeout{** the default language instead.}%
\else
\language=\csname l@#1\endcsname
\fi
#2}}
\providecommand{\BIBdecl}{\relax}
\BIBdecl

\bibitem{ref:subsclust}
L.~Han and X.-L. Liu, ``Convex relaxation algorithm for a structured
  simultaneous low-rank and sparse recovery problem,'' \emph{Journal of the
  Operations Research Society of China}, vol.~3, no.~3, pp. 363--379, 2015.

\bibitem{ref:matrixrigid}
L.~G. Valiant, \emph{Graph-theoretic arguments in low-level complexity}.\hskip
  1em plus 0.5em minus 0.4em\relax Berlin, Heidelberg: Springer Berlin
  Heidelberg, 1977, pp. 162--176.

\bibitem{ref:face2}
N.~S. Aybat, D.~Goldfarb, and S.~Ma, ``Efficient algorithms for robust and
  stable principal component pursuit problems,'' \emph{Computational
  Optimization and Applications}, vol.~58, no.~1, pp. 1--29, 2014.

\bibitem{ref:canPCA}
E.~J. Cand\`{e}s, X.~Li, Y.~Ma, and J.~Wright, ``Robust principal component
  analysis?'' \emph{J. ACM}, vol.~58, no.~3, pp. 11:1--11:37, June 2011.

\bibitem{ref:face1}
J.~Wright, G.~Arvind, R.~Shankar, P.~Yigang, and Y.~Ma, ``Robust principal
  component analysis: Exact recovery of corrupted low-rank matrices via convex
  optimization,'' in \emph{Advances in Neural Information Processing Systems
  22}.\hskip 1em plus 0.5em minus 0.4em\relax Curran Associates, Inc., 2009,
  pp. 2080--2088.

\bibitem{ref:bg1}
K.~Min, Z.~Zhang, J.~Wright, and Y.~Ma, ``Decomposing background topics from
  keywords by principal component pursuit,'' in \emph{Proceedings of the 19th
  ACM International Conference on Information and Knowledge Management}.\hskip
  1em plus 0.5em minus 0.4em\relax New York, NY, USA: ACM, 2010, pp. 269--278.

\bibitem{ref:linkpred}
------, ``Estimation of simultaneously sparse and low rank matrices,'' in
  \emph{29th International conference on machine learning}, 2012, pp.
  1351--1358.

\bibitem{ref:graphclust}
R.~K. Vinayak, S.~Oymak, and B.~Hassibi, ``Sharp performance bounds for graph
  clustering via convex optimization,'' in \emph{2014 IEEE International
  Conference on Acoustics, Speech and Signal Processing (ICASSP)}, May 2014,
  pp. 8297--8301.

\bibitem{ref:pcaspringer}
I.~Jolli, \emph{Principal Component Analysis}.\hskip 1em plus 0.5em minus
  0.4em\relax Springer-Verlag, 1986.

\bibitem{ref:icabook}
A.~Hyvarien, J.~Karhunen, and E.~Oja, \emph{Independent Component
  Analysis}.\hskip 1em plus 0.5em minus 0.4em\relax John Wiley and Sons, New
  York, 2001.

\bibitem{ref:sparseca}
Y.~Li, A.~Cichocki, and S.~Amari, ``Sparse component analysis for blind source
  separation with less sensors than sources,'' in \emph{Proc. Int. Conf.
  Independent Component Analysis (ICA)}, 2003, pp. 89--94.

\bibitem{ref:sparseca2}
B.~A. Olshausen and D.~J. Field, ``Natural image statistics and efficient
  coding,'' \emph{Network Comp. Neural Syst.}, vol.~7, no.~2, pp. 333--339,
  1996.

\bibitem{ref:tesnor}
H.~Rauhut, R.~Schneider, and Z.~Stojanac, ``Low rank tensor recovery via
  iterative hard thresholding,'' in \emph{In Proc. 10th International
  Conference on Sampling Theory and Applications}, 2013.

\bibitem{ref:niht}
T.~Blumensath and M.~E. Davies, ``Normalized iterative hard thresholding:
  Guaranteed stability and performance,'' \emph{IEEE Journal of Selected Topics
  in Signal Processing}, vol.~4, no.~2, pp. 298--309, April 2010.

\bibitem{ref:mtxccom}
J.~Tanner and K.~Wei, ``Normalized iterative hard thresholding for matrix
  completion,'' \emph{SIAM Journal on Scientific Computing}, vol.~35, no.~5,
  pp. S104--S125, 2013.

\bibitem{ref:imat}
F.~Marvasti, A.~Amini, F.~Haddadi, M.~Soltanolkotabi, B.~Khalaj, A.~Aldroubi,
  S.~Sanei, and J.~Chambers, ``A unified approach to sparse signal
  processing,'' \emph{EURASIP Journal on Advances in Signal Processing}, 2012.

\bibitem{ref:inexact}
Z.~Lin, M.~Chen, and Y.~Ma, ``The augmented {L}agrange multiplier method for
  exact recovery of corrupted low-rank matrices,'' \emph{UIUC Technical Report
  UILU-ENG-09-2214, https://arxiv.org/abs/1009.5055v3}, October 2013.

\bibitem{ref:alm1}
X.~Yuan and J.~Yang, ``Sparse and low-rank matrix decomposition via alternating
  direction methods,'' \emph{Optimization Online}, November 2009.

\bibitem{ref:nominalsequence}
Y.~Sheikh and M.~Shah, ``Bayesian modeling of dynamic scenes for object
  detection,'' \emph{IEEE Transactions on Pattern Analysis and Machine
  Intelligence}, vol.~27, no.~11, pp. 1778--1792, November 2005.

\bibitem{ref:exyale}
A.~S. Georghiades, P.~N. Belhumeur, and D.~J. Kriegman, ``From few to many:
  illumination cone models for face recognition under variable lighting and
  pose,'' \emph{IEEE Transactions on Pattern Analysis and Machine
  Intelligence}, vol.~23, no.~6, pp. 643--660, June 2001.

\end{thebibliography}
\balance

\end{document}